\documentclass[11pt]{article}
\usepackage{amsmath}            
\usepackage{amsfonts,amssymb}   

\def\softd{{\leavevmode\setbox1=\hbox{d}%
    \hbox to 1.05\wd1{d\kern-0.2ex{\char039}\hss}}}

\DeclareMathOperator{\alf}{alph}

\newtheorem{theorem}{Theorem}
\newtheorem{cor}{Corollary}

\newtheorem{lemma}{Lemma}

\begin{document}

\title{IDENTITIES IN BRANDT SEMIGROUPS, REVISITED}
\author{ Mikhail V. Volkov\\ Ural Federal University, Ekaterinburg, Russia, m.v.volkov@urfu.ru}
\date{}

\maketitle

\begin{abstract}
We present a new proof for the main claim made in the author's paper ``On the identity bases of Brandt semigroups'' (Ural.\ Gos.\ Univ.\ Mat.\ Zap. \textbf{14}, no.1 (1985), 38--42); this claim provides an identity basis for an arbitrary Brandt semigroup over a group of finite exponent. We also show how to fill a gap in the original proof of the claim in loc. cit.
\end{abstract}

\section{Introduction}
\label{sec:intro}

We assume the reader's acquaintance with the concepts of an identity and an identity basis as well as other rudiments of the theory of varieties; they all may be found, e.g., in~\cite[Chapter~II]{BuSa81}. Our paper deals with identity bases of a certain species of semigroups which we introduce now.

Let $G$ be a group, $I$ a set with at least 2 elements, and $0$ a ``fresh'' symbol that does not belong to $G\cup I$. We define a multiplication on the set $B(G,I)=I\times G\times I\cup\{0\}$ as follows:
\begin{gather*}
(i,g,j)(k,h,\ell)=
\begin{cases}
(i,gh,\ell)&\text{if $j=k$},\\
0&\text{otherwise},
\end{cases}\ \text{for all $i,j,k,\ell\in I$ and all $g,h\in G$,}\\
0x=0,\  x0=0\ \text{ for all $x\in B(G,I)$}.
\end{gather*}
It is easy to verify that the multiplication is associative so that $B(G,I)$ becomes a semigroup. The semigroup is called the \emph{Brandt semigroup over the group $G$}, and the group $G$ in this context is referred to as the \emph{structure group} of $B(G,I)$ while $I$ is called the \emph{index set}.

Recall that an element $a$ of a semigroup $S$ is said to be \emph{regular} if there exists an element $b\in S$ satisfying $aba=a$ and $bab=b$; it is common to say that $b$ is an \emph{inverse of $a$}. A semigroup is called \emph{regular}  [respectively, \emph{inverse}] if every its element has an inverse [respectively, a unique inverse]. The semigroup $B(G,I)$ is inverse: one can easily check that for all $i,j\in I$ and all $g\in G$, the unique inverse of $(i,g,j)$ is $(j,g^{-1},i)$ and the unique inverse of 0 is 0.

Brandt semigroups arose from a concept invented by Brandt~\cite{Brandt27} in his studies on composition of quaternary quadratic forms; a distinguished role played by Brandt semigroups in the structure theory of inverse semigroups was revealed by Clifford~\cite{Clifford:1942} and Munn~\cite{Munn:1957}. From the varietal viewpoint, Brandt semigroups are of importance as well (see, e.g., \cite[Section~7]{ShSu89}), and this justifies the study of their identities. Since Brandt semigroups happen to be inverse, there is a bifurcation in this study: along with plain identities $u=v$, in which the terms $u$ and $v$ are plain semigroup words, that is, products of variables, one can consider also inverse identities whose terms involve both multiplication and the unary operation of taking the inverse. We notice that even though plain identities form a special instance of inverse ones, this does not imply that the study of the former fully reduces to the study of the latter; see Section~\ref{sec:problems} for a more detailed discussion.

Kleiman~\cite{Klei77} comprehensively analyzed inverse identities of Brandt semigroups. In particular, he showed how to derive a basis for such identities of $B(G,I)$ from any given identity basis of the group $G$. Mashevitzky~\cite{Mash79} gave a characterization of the set of all plain identities holding in a given Brandt semigroup modulo the plain identities of its structure group. Trahtman~\cite{Tr81} found a basis for plain identities of the 5-element Brandt semigroup $B_2$ in which the construction $B(G,I)$ results provided that $G$ is the trivial group $E$ and $|I|=2$; this basis consists of the following identities:
\begin{equation}
\label{eq:B2}
x^2=x^3,\ xyx=xyxyx,\ x^2y^2=y^2x^2.
\end{equation}
This fact was frequently cited and used in many applications, including quite important ones such as the positive solution to the finite basis problem for 5-element semigroups~\cite{Tr83,Tr91,Lee13}.

In \cite{Vo85}, the present author applied Kleiman's result from~\cite{Klei77} along with a generalization of Trahtman's argument from~\cite{Tr81} in order to obtain a basis of plain identities for an arbitrary Brandt semigroup over a group of finite exponent. Recall that a group $G$ is said to be of \emph{finite exponent} if there exists a positive integer $n$ such that $g^n=1$ for all $g\in G$. The least number $n$ with this property is called the \emph{exponent} of $G$. Clearly, if $G$ is a group of exponent $n>1$, then $g^{-1}=g^{n-1}$ for all $g\in G$, whence every terms, which is built from certain variables with the help of the unary operation of taking the inverse along with the multiplication, is equal in $G$ to a semigroup word over the same variables. In particular, identities of $G$ (both inverse and plain) admit a basis $\{w_\lambda=1\}_{\lambda\in\Lambda}$ such that each $w_\lambda$ is a plain semigroup word; we refer to such a basis as a \emph{positive identity basis} of $G$. The following is the main result of~\cite{Vo85}:

\begin{theorem}
\label{thm:basis}
Let $G$ be a group of exponent $n>1$, $\{w_\lambda=1\}_{\lambda\in\Lambda}$ a positive identity basis of $G$, and $I$ a set with at least $2$ elements. The identities
\begin{gather}
w_\lambda^2=w_\lambda\quad(\lambda\in\Lambda),\label{eq:idempotent}\\
x^2=x^{n+2},\label{eq:exp n}\\
xyx=(xy)^{n+1}x,\label{eq:exp n_red}\\
x^ny^n=y^nx^n\label{eq:commut}
\end{gather}
constitute a basis for plain identities of the Brandt semigroup $B(G,I)$.
\end{theorem}

This result also has some important consequences, e.g., it implies a classification of finite inverse semigroups whose plain identities admit a finite  basis (\cite[Corollary~3]{Vo85}, see also Section~\ref{sec:problems}).

For more than 25 years there was no doubt in the validity of Trahtman's argument in~\cite{Tr81} until Reilly~\cite{Rei08} observed that the argument in fact contained a lacuna. Nevertheless, the claim made in~\cite{Tr81} turned out to persist since Reilly managed to prove that the identities~\eqref{eq:B2} do form a basis for plain identities of the semigroup $B_2$, see~\cite[Theorem 5.4]{Rei08}. A crucial step in Reilly's proof employs a solution to the word problem in the free objects of the variety generated by $B_2$; this solution (first provided by Mashevitsky in~\cite{Mash79}) has quite a complicated formulation and a somewhat bulky justification. Independently and simultaneously, Lee and the present author~\cite{LeeVo07} invented an alternative way to save Trahtman's claim; their approach bypassed the word problem and resulted in a proof which was short and rather straightforward modulo an elementary yet powerful argument known as Kublanovskii's Lemma, see \cite[Lemma~3.2]{HKMST}. This technique stems from the present author's paper~\cite{Vo05}.

Since the proof of Theorem~\ref{thm:basis} in~\cite{Vo85} uses a version of Trahtman's argument, it suffers from the same problem as the proof in~\cite{Tr81}, and therefore, cannot be considered as truly complete. In fact, the gap in the proof in~\cite{Vo85} can be filled, and we show below how to rescue that proof. However, the main aim of the present paper is to present a new proof of Theorem~\ref{thm:basis}; this new proof follows the approach in~\cite{Vo05,LeeVo07} and relies on a suitable version of  Kublanovskii's Lemma. We have made a fair effort to make our proof self-contained so that, in particular, it should be understandable without any acquaintance with~\cite{Vo85} as a whole nor with specific results therein.

\section{Preliminaries}
\label{sec:prelim}
Here we collect a few auxiliary results that we need; they all either are known or constitute minor variations of known facts. Some of these results and/or their proofs involve certain concepts of semigroup theory, which all can be found in the early chapters of any general semigroup theory text such as, e.g., \cite{ClPr61,How95}.

\begin{lemma}
\label{lem:b2vg}
Let $G$ be an arbitrary group, $I$ a set with at least $2$ elements. An identity $u=v$ holds in the Brandt semigroup $B(G,I)$ if and only if $u=v$ holds in both $G$ and the
$5$-element Brandt semigroup $B_2$.
\end{lemma}

\emph{Proof.}
This was established in \cite[Lemma~5]{Klei77} for inverse identities. As plain identities are special instances of inverse ones, the claim holds for plain identities as well.  \hfill$\square$

\begin{lemma}
\label{lem:idempotent}
Let $G$ be a group and $I$ a set such that $|G|,|I|\ge2$. If $G$ satisfies the identity $w=1$ where $w$ is a semigroup word, then the Brandt semigroup $B(G,I)$ satisfies the identity $w^2=w$.
\end{lemma}

\emph{Proof.}
This fact was also mentioned in \cite[p.~214]{Klei77} for inverse identities, and we could have specialized it to plain identities as we did in the proof of Lemma~\ref{lem:b2vg}. However, the proof in~\cite{Klei77} is only briefly outlined, and the outline involves several advanced notions and results from the theory of inverse semigroups. For the sake of completeness, we provide here a direct and elementary proof.

Clearly, $G$ satisfies the identity $w^2=w$. In view of Lemma~\ref{lem:b2vg} it remains to verify that the identity holds in the semigroup $B_2$. Let $\mathcal{P}(G)$ stand for the set of all non-empty subsets of $G$. We define a multiplication $\cdot$ on the set $\mathcal{P}(G)\times G$ by the following rule: for $A,B\subseteq G$, $g,h\in G$,
\begin{equation}
\label{eq:mult}
(A,g)\cdot(B,h)=(A\cup gB,gh)\ \text{ where }\ gB=\{gb : b\in B\}.
\end{equation}
It is routine to verify that $\cdot$ is associative so that $(\mathcal{P}(G)\times G,\cdot)$ becomes a semigroup which, for brevity, we denote by $S$.

Let $\alf(w)$ denote the set of variables that occur in $w$. If we evaluate the variables $x_1,x_2,\dots\in\alf(w)$ at some elements $(A_1,g_1),(A_2,g_2),\dots$ of $S$ and calculate the corresponding value of $w$, then, according to \eqref{eq:mult}, we get an element of the form $(A,w(g_1,g_2,\dots))$ for a certain set $A\in\mathcal{P}(G)$. Since the identity $w=1$ holds in $G$, we have $w(g_1,g_2,\dots)=1$, so that the value is actually of the form $(A,1)$. Clearly, $(A,1)\cdot(A,1)=(A\cup A,1)=(A,1)$ for every $A\in\mathcal{P}(G)$, whence $S$~satisfies the identity $w^2=w$.

Consider the Brandt semigroup $B(E,G)$ over the trivial group $E=\{1\}$; observe that here we make the set $G$ play the role of the index set! Let $J=\{(A,g)\in S : |A|\ge 2\}$ and define a map $\varphi\colon S\to B(E,G)$, letting $s\varphi=0$ for all $s\in J$ and $(\{a\},g)\varphi=(a,1,g^{-1}a)$ for all $(\{a\},g)\in S\setminus J$. It is easy to see that $\varphi$ is onto: indeed, an arbitrary triple $(k,1,\ell)\in B(E,G)\setminus\{0\}$, where $k,\ell\in G$, has a unique preimage in $S\setminus J$, namely, the pair $(\{k\},k\ell^{-1})$, and for $0$, every element of $J$ is a preimage. Let us verify that $\varphi$ is a semigroup homomorphism. Clearly, $(s\cdot t)\varphi=0=s\varphi\,t\varphi$ whenever at least one of the elements $s$ and $t$ lies in $J$. For $(\{a\},g),(\{b\},h)\in S\setminus J$, we have
\begin{gather*}
\bigl((\{a\},g)\cdot(\{b\},h)\bigr)\varphi=\bigl((\{a,gb\},gh)\bigr)\varphi=\left\{
\begin{array}{ccc}
\text{[if } a=gb]&(a,1,(gh)^{-1}a) &=\\
\text{[if } a\ne gb]&0 &=
\end{array}
\right.\\
\left.\begin{array}{cl}
(a,1,h^{-1}b) &[\text{if }  g^{-1}a=b]\\
0 &[\text{if } g^{-1}a\ne b]
\end{array}
\right\}=(a,1,g^{-1}a)(b,1,h^{-1}b)=(\{a\},g)\varphi\,(\{b\},h)\varphi.
\end{gather*}
Summing up the established properties of $\varphi$, we conclude that the Brandt semigroup $B(E,G)$ is a homomorphic image of the semigroup $S$, and therefore, $B(E,G)$ also  satisfies the identity $w^2=w$.

Since $|G|\ge 2$, we can fix any 2-element subset $K$ in $G$ and ``restrict'' $B(E,G)$ to $K$, that is, consider the subsemigroup $\{(k,1,\ell)\in B(E,G): k,\ell\in K\}\cup\{0\}$ of $B(E,G)$. Then the identity $w^2=w$ holds in this subsemigroup, which clearly is isomorphic to $B_2$.
\hfill$\square$\\[.25ex]%

\emph{Remark~1.} The reader may wonder why Lemma~\ref{lem:idempotent} could not have been proved by a direct evaluation of the word $w$ in the Brandt semigroup $B(G,I)$. The difficulty is that on this way one should have verified that $w$ and $w^2$ take value 0 under the same evaluations of the variables from $\alf(w)$ in $B(G,I)$. Of course, not every word $w$ enjoys this property so that one should have analyzed the structure of $w$, relying entirely on the fact that the identity $w=1$ holds in some non-trivial group. Such an analysis is possible but is rather cumbersome (it amounts to characterizing words $w$ such that the normal closure of $w$ in the free group on the set $\alf(w)$ coincides with the whole group).

\begin{lemma}
\label{lem:splitting}
Let $G$ be a group and $I$ a set  with at least $2$ elements. If the Brandt semigroup $B(G,I)$ satisfies an identity $u=v$ such that $u=u'yu''$ where $y$ is a variable with $y\notin\alf(u'u'')$ and $\alf(u')\cap\alf(u'')=\varnothing$, then $v$ can be decomposed as $v=v'yv''$ with $\alf(v')=\alf(u')$, $\alf(v'')=\alf(u'')$, and the identities $u'=v'$ and $u''=v''$ hold in $B(G,I)$.
\end{lemma}

\emph{Proof.}
One could have deduced Lemma~\ref{lem:splitting} by combining Proposition~3.2(ii) of~\cite{LeeVo07} with its left-right dual. However, since the proof of Proposition~3.2(ii) is omitted in~\cite{LeeVo07}, we prefer to prove the lemma from scratch by a straightforward argument.

Fix two elements $k,\ell\in I$. Suppose that there exists a variable that occurs in only one of the words $u$ and $v$. Evaluating this variable at 0 and other variables at $(k,1,k)$, we get that one of the words $u$ and $v$ takes value 0 while the value of the other is $(k,1,k)$, a contradiction. Hence, $\alf(u)=\alf(v)$. Define an evaluation
$\zeta\colon\alf(u)\to B(G,I)$ as follows:
\[
x\zeta=\begin{cases}
(k,1,k) &\text{if } x\in\alf(u'),\\
(k,1,\ell) &\text{if } x=y,\\
(\ell,1,\ell) &\text{if } x\in\alf(u'').
\end{cases}
\]
Using the multiplication rules of $B(G,I)$, one readily calculates that the value of the word $u$ under $\zeta$ is $(k,1,\ell)$. Since $B(G,I)$ satisfies the identity $u=v$, the value of $v$ under $\zeta$ is $(k,1,\ell)$ as well. This value is a product of the triples $(k,1,k)$, $(k,1,\ell)$, and $(\ell,1,\ell)$ in the same order in which the variables from $\alf(u')$, the variable $y$, and the variables from $\alf(u'')$, respectively, occur in the word $v$. Fix an occurrence of $y$ in $v$ and let $v'y$ be the prefix of $v$ ending with this occurrence and $yv''$ the suffix of $v$ starting with this occurrence. Then $v=v'yv''$. Since
\begin{multline*}
(k,1,\ell)(k,1,\ell)=(k,1,\ell)(k,1,k)=(k,1,k)(\ell,1,\ell)=\\
(\ell,1,\ell)(k,1,\ell)=(\ell,1,\ell)(k,1,k)=0,
\end{multline*}
none of the factors $y^2,yx,xz,zy,zx$ with $x\in\alf(u')$ and $z\in\alf(u'')$ may occur in $v$. Therefore, every variable that appears in $v'$ must come from $\alf(u')$ while every variable that appears in $v''$ must belong to $\alf(u'')$. We see that  $\alf(v')\subseteq\alf(u')$, $\alf(v'')\subseteq\alf(u'')$, and from the equality  $\alf(u)=\alf(v)$ shown above, we conclude that  $\alf(v')=\alf(u')$, $\alf(v'')=\alf(u'')$.

It remains to verify that the identities $u'=v'$ and $u''=v''$ hold in $B(G,I)$. The semigroup $B(G,I)$ is inverse, and every inverse semigroup is isomorphic to its left-right dual via the bijection that maps each element to its unique inverse. Therefore $B(G,I)$ satisfies an identity $p=q$ if and only if it satisfies its mirror image $\overleftarrow{p}=\overleftarrow{q}$, where $\overleftarrow{w}$ denotes the word $w$ read backwards. In view of this symmetry, it suffices to verify that $u'=v'$ holds in $B(G,I)$. Arguing by contradiction, consider an evaluation $\varphi\colon\alf(u')\to B(G,I)$ such that the values of $u'$ and $v'$ under $\varphi$ are different. Then one of these values is not equal to 0; assume, for certainty, that the value of $u'$ is some triple $(i,g,j)\in B(G,I)\setminus\{0\}$. We extend $\varphi$ to an evaluation $\psi\colon\alf(u)\to B(G,I)$, letting $x\psi=x\varphi$ for all $x\in\alf(u')$ and $y\psi=z\psi=(j,1,j)$ for all $z\in\alf(u'')$. The value of $u$ under $\psi$ is $(i,g,j)(j,1,j)=(i,g,j)$; we aim to show that the value of $v$ under $\psi$ is different from $(i,g,j)$. Indeed, if the value of $v'$ under $\varphi$ is 0, so is the value of $v$ under $\psi$. If the value of $v'$ under $\varphi$ is a triple $(i',g',j')\ne(i,g,j)$, then the value of $v$ under $\psi$ is
\[
(i',g',j')(j,1,j)=\begin{cases}
(i',g',j)&\text{if $j'=j$},\\
0&\text{if $j'\ne j$},
\end{cases}\ \ne\ (i,g,j).
\]
This contradicts the premise of $u=v$ holding in $B(G,I)$.\hfill$\square$\\[.25ex]%

A $[0]$-\emph{minimal ideal} of a semigroup $S$ is its minimal (with respect to the set inclusion) non-zero ideal if $S$ has a zero and its least ideal otherwise.
A non-trivial semigroup $S$ is $[0]$-\emph{simple} if $S=S^2$ and $S$ is a $[0]$-minimal ideal of itself. A $[0]$-simple semigroup is \emph{completely} $[0]$-\emph{simple} if it contains an idempotent $e$ such that every idempotent $f$ satisfying $ef=fe=f$ is equal to either $e$ or $0$.

\begin{lemma}
\label{lem:minimal}
If a semigroup satisfies the identities \eqref{eq:exp n_red} and \eqref{eq:commut} for some $n\ge 1$, then every its $[0]$-minimal ideal that contains a regular element is
an inverse completely $[0]$-simple semigroup.
\end{lemma}

\emph{Proof.}
It suffices to combine a few standard facts of semigroup theory. First, in any semigroup, a $[0]$-minimal ideal with a regular element is a $[0]$-simple semigroup, see \cite[Theorem~2.29]{ClPr61} or \cite[Proposition~3.1.3]{How95}. Second, every [0]-simple semigroup that satisfies \eqref{eq:exp n_red} is completely $[0]$-simple; this is a special case of Munn's theorem, see \cite[Theorem~2.55]{ClPr61} or \cite[Theorem~3.2.11]{How95}. Each completely $[0]$-simple semigroup is regular, and a regular semigroup with commuting idempotents is inverse, see \cite[Theorem~1.17]{ClPr61} or \cite[Theorem~5.1.1]{How95}. It remains to observe that idempotents commute in every semigroup satisfying \eqref{eq:commut}.
\hfill$\square$\\[.25ex]%

We say that a map $\varphi\colon S\to T$ \emph{separates elements $a,b\in S$} if $a\varphi\ne b\varphi$.

\begin{lemma}
\label{lem:kublanovskii}
If  a semigroup $S$ satisfies the identities \eqref{eq:exp n_red} and \eqref{eq:commut} for some $n\ge 1$, then any distinct regular elements $a,b\in S$ are separated by a homomorphism of $S$ onto an inverse completely $[0]$-simple semigroup.
\end{lemma}

\emph{Proof.}
This is a version of Kublanovskii's Lemma \cite[Lemma~3.2]{HKMST} adapted for the purposes of the present paper.  For the reader's convenience, we provide a complete proof, even though it quite closely follows the proof of Kublanovskii's Lemma in~\cite{HKMST}.

For each regular element $z\in S$, we let $I_z=\left\{u\in S : z\notin SuS\right\}$. Observe that $z\notin I_z$: indeed, if $z'$ is an inverse of $z$, we have $z=zz'zz'z\in SzS$. The set $I_z$ may be empty but if it is not empty, it forms an ideal of $S$. Indeed, $SutS\subseteq SuS$ and $StuS\subseteq SuS$ for any $u,t\in S$, and hence, if $u$ lies in $I_z$, so do $ut$ and $tu$ for every $t\in S$. Define the following equivalence relation on $S$:
\[
 x\equiv y\!\!\pmod{I_z}\ \text{ if and only if either } x=y \text{ or } x,y\in I_z.
\]
Clearly, it is just the equality relation if $I_z$ is empty; otherwise it is nothing but the Rees congruence $\iota_z$ corresponding to the ideal $I_z$. Now define a further
equivalence relation $\rho_z$ on $S$ as follows:
\[
\rho_z =\left\{\left( x,y\right)\in S\times S :  xt\equiv yt\!\!\pmod{I_z}\ \text{ for every } t\in SzS\right\}.
\]
It can be easily verified that $\rho _z$ is a congruence on $S$; in fact, as observed in~\cite{HKMST}, $\rho _z$ is the kernel of the so-called Sch\"utzenberger representation for $S$, see \cite[Section~3.5]{ClPr61}.

Clearly, $\rho_z=S\times S$ if $z=0$. Now we aim to prove the following claim: \emph{if $z\ne 0$, then the quotient $S/\rho_z$ is an inverse completely $[0]$-simple semigroup}.

If $I_z\ne\varnothing$, the congruence $\rho_z$ contains the Rees congruence $\iota_z$. Then we may substitute $S$ by its quotient $S/\iota_z$ as the quotient also satisfies the identities \eqref{eq:exp n_red} and \eqref{eq:commut}; in other words, we may (and will) assume that either $I_z=\varnothing$ or $I_z=\{0\}$. Then by the definition of the set $I_z$, every non-zero element $u\in SzS$ must fulfil $z\in SuS$ whence $SuS=SzS$. We see that $SzS$ is a $[0]$-minimal ideal of $S$; as $SzS$ contains $z$ which is a regular element, Lemma~\ref{lem:minimal} applies showing that $SzS$ is an inverse completely $[0]$-simple semigroup. So is any homomorphic image of $SzS$; in particular, so is the image of $SzS$ in the quotient semigroup $S/\rho_z$.  Therefore, it remains to show that the image of $S$ in $S/\rho_z$ coincides with that of $SzS$, which means that for each $x\in S$, there exists $y\in SzS$ such that $\left(x,y\right)\in\rho_z$.

If $x\in SzS$, there is nothing to prove. If $x\notin SzS$, then in particular, $x\notin I_z$ whence $z=pxq$ for some $p,q\in S$. We have $z=pxqz'pxq$, where, as above, $z'$ stands for an inverse of $z$. Put $w=qz'p$; then $w\in SzS$ because $z'=z'zz'\in SzS$ and $xwx\ne0$ because $z=pxwxq\ne0$. Now take an arbitrary element $t\in SzS$. We have already noticed (in the preceding paragraph) that $SuS=SzS$ for   every non-zero element $u\in SzS$. Applying this to $u=xwx$, we conclude that $t=rxwxs$ for some $r,s\in S$. Now we have the following chain of equalities:
\begin{align*}
xt=xrxwxs&=(xr)^{n+1}(xw)^{n+1}xs&&\text{by applying \eqref{eq:exp n_red} to $xrx$ and $xwx$}\\
         &=xr(xr)^n(xw)^nxwxs&&\\
         &=xr(xw)^n(xr)^nxwxs&&\text{by applying \eqref{eq:commut}}\\
         &=xr(xw)^n(xr)^{n-1}xrxwxs&&\\
         &=xr(xw)^n(xr)^{n-1}xt.&&
\end{align*}
We see that $\left(x,xr(xw)^n(xr)^{n-1}x\right)\in\rho _z$, and the element $xr(xw)^n(xr)^{n-1}x$ lies in the ideal $SzS$ because so does $w$. Thus, $xr(xw)^n(xr)^{n-1}x$ can play the role of $y$, and our claim is proved.

Now we are ready to complete the proof of the lemma. Given an arbitrary pair $(a,b)$ of distinct regular elements is $S$, we will show that at least one of the congruences $\rho_a$ and $\rho_b$ excludes $(a,b)$. Then the natural homomorphism of $S$ onto the quotient over this congruence separates $a$ and $b$, and the quotient is an inverse completely $[0]$-simple semigroup by the claim just proved. (One has to take into account that if a congruence of the form $\rho_z$ excludes some pair, then $z\ne0$ and the claim applies.)

If $a\notin SbS$, then $b\in I_a$. Let $a'$ be an inverse of $a$. We have then $a'a\in SaS$ and $a(a'a)=a\notin I_a$ while $b(a'a)\in I_a$ since $I_a$ is an ideal. Hence $(a,b)\notin\rho_a$. Similarly, if $b\notin SaS$, we have $(a,b)\notin\rho_b$. Now suppose that $a\in SbS$ and $b\in SaS$. In this case, $SaS=SbS$ and $a,b\notin I_a=I_b$. If we assume that $(a,b)\in\rho_a$, then for every element $t\in SaS$ such that either $at\notin I_a$ or $bt\notin I_a$, we must have $at=bt$. In particular, the latter equality must hold for $t=a'a$ since $a(a'a)=a\notin I_a$ and for $t=b'b$, where $b'$ is an inverse of $b$, since $b(b'b)=b\notin I_a$. Taking into account that both $a'a$ and $b'b$ are idempotents and that idempotents commute in every semigroup satisfying the identity \eqref{eq:commut}, we have
\begin{multline*}
a=a(a'a)=b(a'a)=b(b'b)(a'a)=a(b'b)(a'a)=\\
a(a'a)(b'b)=a(b'b)=b(b'b)=b,
\end{multline*}
a contradiction.\hfill$\square$\\[.25ex]%

\emph{Remark~2.} One can call our Lemma~\ref{lem:kublanovskii} ``Kublanovskii's Lemma with commuting idempotents''. The presence of the identity \eqref{eq:commut} ensures that idempotents commute, and this streamlines the proof. The most important simplification in comparison with the proof of Kublanovskii's Lemma in~\cite{HKMST} is that we manage to avoid invoking, along  with the congruences $\rho_a$ and $\rho_b$, their dual versions, that is, the kernels of the corresponding Sch\"utzenberger anti-representations.\\[.25ex]

If $S$ is an arbitrary semigroup and $0$ is a ``fresh'' symbol that does not belong to $S$, we let $S^0$ stand for the semigroup on the set $S\cup\{0\}$ with multiplication that extends the multiplication of $S$ and makes all products involving $0$ be equal to $0$. If $G$ is a group, $G^0$ is known under the (standard though somewhat oxymoronic) name ``\emph{group with zero}''. The following fact is a classical result of semigroup theory, see  \cite[Theorem~3.9]{ClPr61} or \cite[Theorem 5.1.8]{How95}.

\begin{lemma}
\label{lem:inv_c0s}
 An inverse completely $[0]$-simple semigroup is either a group, or a group with zero, or a Brandt semigroup.
\end{lemma}

\section{Proof of Theorem~\ref{thm:basis}}
\label{sec:proof}

Recall that we aim to prove that for every group $G$ of exponent $n>1$ and every set $I$ with at least $2$ elements, the identities \eqref{eq:idempotent}--\eqref{eq:commut} constitute a basis of the plain identities of the Brandt semigroup $B(G,I)$, provided that the set $\{w_\lambda=1\}_{\lambda\in\Lambda}$ is a positive identity basis of $G$.

To start with, observe that the identities \eqref{eq:idempotent}--\eqref{eq:commut} hold in $B(G,I)$. For \eqref{eq:idempotent} this follows from Lemma~\ref{lem:idempotent}. As for the identities \eqref{eq:exp n}--\eqref{eq:commut}, it is obvious that they hold in each group of exponent~$n$. On the other hand, comparing these identities with the identity basis \eqref{eq:B2} of the semigroup $B_2$, one readily sees that they hold in $B_2$ as well. Now the ``if'' part of Lemma~\ref{lem:b2vg} ensures that \eqref{eq:exp n}--\eqref{eq:commut} hold in $B(G,I)$.

Let $\mathbf{A}$ be the semigroup variety defined by the identities \eqref{eq:idempotent}--\eqref{eq:commut} and $\mathbf{B}$ the variety generated by the Brandt semigroup $B(G,I)$. The fact established in the preceding paragraph is equivalent to the inclusion $\mathbf{B}\subseteq\mathbf{A}$ and the theorem being proved means the equality $\mathbf{B}=\mathbf{A}$. Arguing by contradiction, assume that the inclusion is strict. Then there exists an identity that holds in the semigroup $B(G,I)$ but fails in the variety $\mathbf{A}$. We choose an identity $u=v$ with this property and with the least value of $|\alf(u)|$. We first check that the words $u$ and $v$ are repeated, where a word $w$ is called \emph{repeated} if each variable from $\alf(w)$ occurs in a factor of $w$ of the form $ypy$ where $y$ is a variable and $p$ is a (possibly empty) word\footnote{The term ``repeated'' comes from \cite{Tr81,Vo85}; in \cite{LeeVo07} words with this property were called ``semiconnected''.}. It is convenient to have a short name for such factors; let us refer to them as to \emph{cells}.

Assume for a moment that, say, $u$ is not repeated. This means that there exists a variable $y$ that occurs in $u$ but does not occur in any cell of $u$. In particular, $y$ occurs in $u$ exactly once, and moreover, $u=u'yu''$ with $\alf(u')\cap\alf(u'')=\varnothing$. We are in a position to employ Lemma~\ref{lem:splitting} to conclude that $v$ decomposes as $v=v'yv''$ where $\alf(v')=\alf(u')$, $\alf(v'')=\alf(u'')$ and both the identities $u'=v'$ and $u''=v''$ hold in $B(G,I)$. Since $|\alf(u')|,|\alf(u'')|<|\alf(u)|$, our choice of the identity $u=v$ ensures that the identities $u'=v'$ and $u''=v''$ hold in the variety $\mathbf{A}$. However, together they imply the identity $u=v$ that cannot hold in $\mathbf{A}$, a contradiction.

Let $F$ stand for the free semigroup of countable rank and let $\alpha$ denote the fully invariant congruence on $F$ that corresponds to the variety $\mathbf{A}$. Then the quotient semigroup $F/\alpha$ satisfies the identities \eqref{eq:idempotent}--\eqref{eq:commut} and the $\alpha$-classes $u^\alpha=\{w: (w,u)\in\alpha\}$ and $v^\alpha=\{w: (w,v)\in\alpha\}$ are  different in $F/\alpha$. For the next step of our proof we need the following fact:
\begin{lemma}
\label{lem:regular}
Every $\alpha$-class that contains a repeated word is a regular element of $F/\alpha$.
\end{lemma}
We proceed with proving Theorem~\ref{thm:basis} modulo Lemma~\ref{lem:regular} and prove the lemma afterwards.

By Lemma~\ref{lem:regular}, the $\alpha$-classes $u^\alpha$ and $v^\alpha$ are regular elements of $F/\alpha$. Applying Lemma~\ref{lem:kublanovskii}, we conclude that $u^\alpha$ and $v^\alpha$ are separated by an onto homomorphism $\chi\colon F/\alpha\to T$, where $T$ is an inverse completely $[0]$-simple semigroup. Lemma~\ref{lem:inv_c0s} implies the existence of a group $Q$ such that either 1) $T=Q$, or 2) $T=Q^0$, or 3) $T=B(Q,J)$ for some set $J$ with $|J|\ge2$. In any case, $Q$ is a subgroup of a homomorphic image of $F/\alpha$, whence the identities~\eqref{eq:idempotent} hold in $Q$. Clearly, if for some word $w$, a group satisfies the identity $w^2=w$, then the group satisfies the identity $w=1$ as well. Therefore the group $Q$ satisfies the identities $w_\lambda=1$ for all $\lambda\in\Lambda$. Since these identities form a basis for the identities of the structure group $G$ of our semigroup $B(G,I)$, the group $Q$ belongs to the semigroup variety generated by $G$, and hence, to the variety $\mathbf{B}$ generated by $B(G,I)$. The 5-element Brandt semigroup $B_2$ also belongs to $\mathbf{B}$; this follows, for instance from the ``only if'' part of Lemma~\ref{lem:b2vg}. Applying the ``if'' part of Lemma~\ref{lem:b2vg}, we conclude that the Brandt semigroup $B(Q,J)$ lies in $\mathbf{B}$. From this, we have $T\in\mathbf{B}$ as $T$ is isomorphic to a subsemigroup in $B(Q,J)$ in the cases 1) or 2) and $T=B(Q,J)$ in the case 3). In particular, $T$ satisfies the identity $u=v$. However, the composition of the natural homomorphism $F\to F/\alpha$ with the homomorphism $\chi\colon F/\alpha\to T$ gives rise to an evaluation under which the values of the words $u$ and $v$ are different. This contradiction completes the proof of Theorem~\ref{thm:basis} modulo Lemma~\ref{lem:regular}.

\medskip

\emph{Proof} of Lemma~\ref{lem:regular}.
Take any  $\alpha$-class $h$ that contains a repeated word, say, $w$. If some variable $y$ occurs in $w$ only once, then by the definition of a repeated word, $y$ occurs in some cell $zpz$ of $w$, where $p$ is non-empty. Using the identity~\eqref{eq:exp n_red}, we substitute the factor $zpz$ by the factor $(zp)^{n+1}z$ and get a new word in the same $\alpha$-class $h$ in which $y$ occurs at least twice. If this new word still contains some variable $x$ with a single occurrence, we apply the same transformation again, etc. Thus, we may assume that $h$ contains a word $q$ in which every variable occurs at least twice. Now we prove that $h$ contains also a word which is a product of cells, that is, has the form
\begin{equation}
\label{eq:cycle}
y_1p_1y_1\cdot y_2p_2y_2\cdot\ldots\cdot y_kp_ky_k,
\end{equation}
where $y_1,y_2,\dots,y_k$ are variables and $p_1,p_2,\dots,p_k$ are (possibly empty) words. For this, we employ a sort of greedy algorithm. Let $y_1$ be the leftmost variable of the word $q$. If $q$ ends with $y_1$, the word $q$ itself is a cell. Otherwise we find the rightmost occurrence of $y_1$ in $q$ so that $q=y_1p_1y_1\cdot q_1$ where $q_1$ is a non-empty word in which $y_1$ does not occur, and so $|\alf(q_1)|<|\alf(q)|$. Let $y_2$ be the leftmost variable of $q_1$. There are two cases to consider, depending on whether $y_2$ occurs in $q_1$ at least twice or only once. In the former case, we find the rightmost occurrence of $y_2$ in $q_1$ and represent $q$ as $q=y_1p_1y_1\cdot y_2p_2y_2\cdot q_2$, where $y_1,y_2$ do not occur in $q_2$, and so $|\alf(q_2)|<|\alf(q_1)|$. Let us show that $h$ contains a word with a similar structure also in the latter case. Indeed, the variable $y_2$ occurs in $q$ at least twice and if it occurs in $q_1$ only once, then it must occur in $p_1$. Hence, $p_1=ry_2s$ for some (possibly empty) words $r$ and $s$. Then $q$ contains the word $y_2sy_1y_2$ as a factor. Using the identity~\eqref{eq:exp n_red}, we substitute this factor by $(y_2sy_1)^{n+1}y_2$ and transform $q$ into a new word $q'$ in the same $\alpha$-class $h$; this new word can be represented as $q'=y_1p'_1y_1\cdot y_2p'_2y_2\cdot q'_2$, where  $p'_1=r(y_2sy_1)^{n-1}y_2s$, $p'_2=sy_1$, and $q'_2$ is obtained from $q_1$ by removing its leftmost variable. Then $y_1,y_2$ do not occur in $q'_2$, whence $|\alf(q'_2)|<|\alf(q_1)|$. Now we can apply the same procedure to the leftmost variable of $q_2$ or $q'_2$, and so on. On the $i$-th step of the procedure we create a new cell $y_ip_iy_i$ while the yet unprocessed ``remainder'' omits the variables $y_1,\dots,y_i$. Clearly, the procedure terminates after a finite number of steps and yields a word of the form~\eqref{eq:cycle} in the $\alpha$-class $h$.

Now let $h^*$  be the $\alpha$-class that contains the word
\[
(p_ky_k)^{2n-2}p_k\cdot(p_{k-1}y_{k-1})^{2n-2}p_{k-1}\cdot\ldots\cdot(p_1y_1)^{2n-2}p_1.
\]
We show that $h^*$ is an inverse of $h$ by induction on $k$. If $k=1$, that is, $h=(y_1p_1y_1)^\alpha$, the $\alpha$-class $hh^*h$ contains the word
\[
y_1p_1y_1\cdot(p_1y_1)^{2n-2}p_1\cdot y_1p_1y_1=(y_1p_1)^{2n+1}y_1.
\]
Applying the identity \eqref{eq:exp n} if the word $p_1$ is empty and the identity \eqref{eq:exp n_red} otherwise, we can transform this word to the word $y_1p_1y_1\in h$. Thus, $hh^*h=h$. Similarly, the $\alpha$-class $h^*hh^*$ contains the word
\[
(p_1y_1)^{2n-2}p_1\cdot y_1p_1y_1\cdot(p_1y_1)^{2n-2}p_1=(p_1y_1)^{4n-2}p_1
\]
that can be transformed to $(p_1y_1)^{2n-2}p_1\in h^*$. Hence, $h^*hh^*=h^*$ and thus, $h^*$ is an inverse of $h$.

For the induction step, suppose that $k>1$ and let $f$ and $g$ be the $\alpha$-classes containing the words $y_1p_1y_1$ and $y_2p_2y_2\cdot\ldots\cdot y_kp_ky_k$ respectively. Then $h=fg$, $h^*=g^*f^*$ and, by the induction assumption, $f^*$ and $g^*$ are inverses of $f$ and $g$, respectively. The equalities $ff^*f=f$ and $gg^*g=g$ imply that the $\alpha$-classes $f^*f$ and $gg^*$ are idempotents. Taking into account that the idempotents of $F/\alpha$ commute due to the identity \eqref{eq:commut}, we obtain
\[
\begin{tabular}{p{6cm}cp{6cm}}
\parbox{6cm}{\begin{align*}
hh^*h &=fg\cdot g^*f^*\cdot fg\\
     &=f(gg^*)(f^*f)g\\
     &=f(f^*f)(gg^*)\\
     &=ff^*f\cdot gg^*g\\
     &=fg=h,
\end{align*}}
&\qquad&
\parbox{6cm}{\begin{align*}
h^*hh^* &= g^*f^*\cdot fg\cdot g^*f^*\\
     &=g^*(f^*f)(gg^*)f^*\\
     &=g^*(gg^*)(f^*f)f^*\\
     &=g^*gg^*\cdot f^*ff^*&&\\
     &=g^*f^*=h^*.
\end{align*}}
\end{tabular}
\]
We see that $h^*$ is an inverse of $h$, and the lemma is proved.
\hfill$\square$\\[.25ex]%

Now we are in a position to discuss a gap in the original proof of Theorem~\ref{thm:basis} in~\cite{Vo85} and to explain how the gap can be filled.

The proof of Theorem~\ref{thm:basis} in~\cite{Vo85} develops as follows. As above, it works with $F$, the free semigroup of countable rank, and $\alpha$, the fully invariant congruence on $F$ that corresponds to the variety $\mathbf{A}$ defined by the identities~\eqref{eq:idempotent}--\eqref{eq:commut}. In the quotient semigroup $F/\alpha$, one considers the set $H$ of all $\alpha$-classes containing a repeated word. Obviously, the product of two repeated words is a repeated word whence $H$ is a subsemigroup of $F/\alpha$. The idempotents of $H$ commute because $H$, being a subsemigroup of $F/\alpha$, satisfies the identity\eqref{eq:commut}. By Lemma~\ref{lem:regular} (which appears in~\cite{Vo85} as a part of the proof of Theorem~\ref{thm:basis}), $H$ is regular. Now one can apply the textbook fact that a regular semigroup with commuting idempotents is inverse, see \cite[Theorem~1.17]{ClPr61} or \cite[Theorem~5.1.1]{How95}. Thus, $H$ is an inverse subsemigroup of $F/\alpha$. At this point, the proof under discussion invokes the main result from Kleiman's paper~\cite{Klei77}, which implies that the identities~\eqref{eq:idempotent}--\eqref{eq:commut} form a basis for the inverse identities of the Brandt semigroup $B(G,I)$. In particular, these identities hold in $B(G,I)$ whence $\mathbf{A}\supseteq\mathbf{B}$, where as above, $\mathbf{B}$ stands for the variety generated by $B(G,I)$. In the language of fully invariant congruences this means that $\alpha\subseteq\beta$, where $\beta$ denotes the fully invariant congruence on $F$ that corresponds to the variety $\mathbf{B}$. Let $\beta/\alpha$ be the induced congruence on $F/\alpha$ so that $\left(F/\alpha\right)/\left(\beta/\alpha\right)\cong F/\beta$. The rest of the proof relies on the following claim: \emph{the congruence $\beta/\alpha$ separates the elements of the subsemigroup $H$}, that is, $\beta/\alpha$ restricted to $H$ is the equality relation. In~\cite{Vo85} this claim is justified by observing that $H$ lies in the variety $\mathbf{B}$---this follows from the fact that $H$ is inverse and satisfies the identities~\eqref{eq:idempotent}--\eqref{eq:commut} which, according to the quoted result from~\cite{Klei77}, define the variety of inverse semigroups generated by $B(G,I)$. However, the justification is not sufficient. The membership $H\in\mathbf{B}$ only guarantees that the \textbf{least} element in the set $\Gamma$ of all congruences $\gamma$ on $H$ with $H/\gamma\in\mathbf{B}$ is the equality relation; while $\beta/\alpha$ restricted to $H$ is a congruence in $\Gamma$, it is not immediately clear that the restriction is indeed the least element in $\Gamma$.

Let us show that the italicized claim does hold. Arguing by contradiction, assume that some distinct elements $a,b\in H$ satisfy $(a,b)\in\beta/\alpha$. Since $a$ and $b$ are distinct regular elements of the semigroup $F/\alpha$, which satisfies the identities \eqref{eq:exp n_red} and \eqref{eq:commut}, Lemma~\ref{lem:kublanovskii} applies. Thus, $a$ and $b$ are separated by an onto homomorphism $\chi\colon F/\alpha\to T$, where $T$ is an inverse completely $[0]$-simple semigroup. Arguing as in the last paragraph of the above proof of Theorem~\ref{thm:basis} modulo Lemma~\ref{lem:regular}, one can show that $T$ lies in the variety $\mathbf{B}$. Then the homomorphism $\chi$ must factor through the natural homomorphism $\eta\colon F/\alpha\to F/\beta$ because $F/\beta$ is the $\mathbf{B}$-free semigroup of countable rank. However, $a\eta=b\eta$ since $(a,b)\in\beta/\alpha$ while $a\chi\ne b\chi$, a contradiction.

\section{Corollaries and discussions}
\label{sec:problems}

For the reader's convenience, we reproduce the main corollaries of Theorem~\ref{thm:basis}, following~\cite{Vo85}.  The first of them specializes Theorem~\ref{thm:basis}, providing an explicit identity basis for Brandt semigroups over abelian groups of finite exponent.

\begin{cor}[{\mdseries\cite[Corollary 1]{Vo85}}]
\label{cor:abelian}
Let $G$ be an abelian group of exponent\linebreak $n>1$ and $I$ a set with at least $2$ elements. The identities \eqref{eq:exp n}, \eqref{eq:exp n_red}, and
\begin{gather}
\label{eq:abelian1}
x^2y^2=y^2x^2,\\
\label{eq:abelian2}
xyxzx=xzxyx
\end{gather}
constitute a basis for plain identities of the Brandt semigroup $B(G,I)$.
\end{cor}

This is in fact a consequence of the proof of Theorem~\ref{thm:basis} rather than the theorem itself. The corresponding arguments were omitted in~\cite{Vo85}; therefore, we provide a proof outline here.

\emph{Proof} (outline).  First, we show that the identities \eqref{eq:exp n}, \eqref{eq:exp n_red}, \eqref{eq:abelian1}, \eqref{eq:abelian2} hold in $B(G,I)$. By the ``if'' part of Lemma~\ref{lem:b2vg}, it suffices to verify that they hold in both $G$ and the 5-element Brandt semigroup $B_2$. Obviously, the identities \eqref{eq:exp n} and \eqref{eq:exp n_red} hold in every group of exponent $n$ while the identities \eqref{eq:abelian1} and \eqref{eq:abelian2} hold in every abelian group. Thus, \eqref{eq:exp n}, \eqref{eq:exp n_red}, \eqref{eq:abelian1}, \eqref{eq:abelian2} hold in $G$. Inspecting the identity basis \eqref{eq:B2}, one readily sees that \eqref{eq:exp n}, \eqref{eq:exp n_red}, \eqref{eq:abelian1} hold in $B_2$. The identity~\eqref{eq:abelian2} also holds in $B_2$ as the following calculation shows:
\begin{align*}
xyxzx&=(xy)^2(xz)^2x&&\text{in view of $xyx=xyxyx$}\\
     &=(xz)^2(xy)^2x&&\text{in view of $x^2y^2=y^2x^2$}\\
     &=xzxyx&&\text{in view of $xyx=xyxyx$.}
\end{align*}

Now we proceed exactly as in the proof of Theorem~\ref{thm:basis}. Denote by $\mathbf{A}$ the semigroup variety defined by the identities \eqref{eq:exp n}, \eqref{eq:exp n_red}, \eqref{eq:abelian1}, \eqref{eq:abelian2} and let $\mathbf{B}$ be the variety generated by the semigroup $B(G,I)$. The fact that $B(G,I)$ satisfies \eqref{eq:exp n}, \eqref{eq:exp n_red}, \eqref{eq:abelian1}, \eqref{eq:abelian2} implies that $\mathbf{B}\subseteq\mathbf{A}$. Assuming that the inclusion is strict, choose an identity $u=v$  with the least value of $|\alf(u)|$ such that $u=v$ holds in $B(G,I)$ but fails in $\mathbf{A}$. Then the words $u$ and $v$ are repeated due to the argument in the 4th paragraph of Section~\ref{sec:proof}.

Let $F$ be the free semigroup of countable rank and $\alpha$ its fully invariant congruence corresponding to the variety $\mathbf{A}$. The $\alpha$-classes $u^\alpha$ and $v^\alpha$ are distinct elements of $F/\alpha$ and, by Lemma~\ref{lem:regular}, they are regular. Then Lemmas~\ref{lem:kublanovskii} and~\ref{lem:inv_c0s} imply that $u^\alpha$ and $v^\alpha$ are separated by an onto homomorphism $\chi\colon F/\alpha\to T$, where $T$ is either a group, or a group with zero, or a Brandt semigroup. Let $Q$ stand for the structure group of $T$ in the latter case and for $T$ or $T\setminus\{0\}$ in the two former cases. Then $Q$ is a subgroup of a homomorphic image of $F/\alpha$, whence the identities~\eqref{eq:exp n} and \eqref{eq:abelian2} hold in $Q$. Clearly, the exponent of every group satisfying~\eqref{eq:exp n} divides $n$ and every group satisfying~\eqref{eq:abelian2} is abelian. Thus, $Q$ is an abelian group of exponent dividing $n$. A well known classification of abelian group varieties (cf. \cite[Theorem 19.5]{Neu37} or \cite[Item 13.51]{Neu67}) ensures that the variety of all abelian groups of exponent dividing $n$ is generated by any abelian group of exponent $n$, in particular, by the structure group $G$ of $B(G,I)$. Thus, $Q$  belongs to the variety generated by $G$, and hence, to the variety $\mathbf{B}$. As the 5-element Brandt semigroup $B_2$ also belongs to $\mathbf{B}$, the ``if'' part of Lemma~\ref{lem:b2vg} implies that every Brandt semigroup over $Q$ lies in $\mathbf{B}$. From this, we have $T\in\mathbf{B}$ whence $T$ must satisfy $u=v$. On the other hand, the composition of the natural homomorphism $F\to F/\alpha$ with the homomorphism $\chi\colon F/\alpha\to T$ separates $u$ and $v$ in $T$, a contradiction.\hfill$\square$\\[.25ex]%

\emph{Remark~3.} We do not know any basis for plain identities of the Brandt semigroup over the infinite cyclic group $\mathbb{Z}$ (or any other abelian group of infinite exponent); moreover, it is not known whether or not the plain identities of this semigroup admit a finite basis. A finite basis for inverse identities of the Brandt semigroup over $\mathbb{Z}$ can be found in~\cite[Corollary~6]{Klei77} or~\cite[Theorem~XII.5.4(iii)]{Petrich84}.\\[.25ex]

In connection with Remark~3, it appears appropriate to discuss in more detail how the \emph{finite basis property}, i.e., the property of a Brandt semigroup $B(G,I)$ to have a finite identity basis, may depend on the type of identities---inverse or plain---under consideration. It turns out that the picture is rather non-trivial here. On the one hand, the additional operation increases the expressivity of the equational language so that the inverse identities of $B(G,I)$ are ``richer'' than the plain ones. This indicates that $B(G,I)$ may have more chances to possess no finite basis for its inverse identities. On the other hand, the inference power of the
language increases too. Hence one can encounter the situation when some identity of $B(G,I)$ does not follow from an identity system $\Sigma$ as a ``plain'' identity but follows from $\Sigma$ as an ``inverse'' identity. This indicates that the inverse identities of $B(G,I)$ may admit a finite basis even if its plain identities do not. The cumulative effect of the trade-off between increased expressivity and increased inference power is hard to predict in general, as the following examples demonstrate\footnote{Our examples are adaptations of known ones (see, e.g., \cite[Section~2]{Vo01}) to the case of Brandt semigroups.}.\\[.25ex]

\emph{Example~1.} Let $G$ be the wreath product of the countably generated free group of exponent 4 with the countably generated free abelian group and $I$ a set with at least $2$ elements. The Brandt semigroup $B(G,I)$ satisfies only trivial plain identities but its inverse identities have no finite basis.\\[.25ex]

\emph{Proof.} The fact that $B(G,I)$ satisfies only trivial plain identities follows from the observation that $G$ contains the countably generated free semigroup as a subsemigroup, see, e.g., \cite{BeSe81}. If we assume that the  inverse identities of $B(G,I)$ admit a finite basis, then appending the identity $xx^{-1}=yy^{-1}$ to the basis would yield a finite basis of group identities of the group $G$. However, by~\cite[Corollary 22.22]{Neu37} $G$ generates the varietal product of the variety of all groups of exponent dividing 4 with the variety of all abelian groups, and by~\cite[Remark~2]{Kl73} this product possesses no finite identity basis, a contradiction. \hfill$\square$\\[.25ex]

In Example~1, an increase in the expressivity of the equational language dominates; now we exhibit an ``opposite'' example in which one sees the effect of an increase in the inference power.\\[.25ex]

\emph{Example~2.} Let $G$ be the direct product of the infinite cyclic group $\mathbb{Z}$ with the group $\mathbb{S}_3$ of all permutations of a $3$-element set and $I$ a set with at least $2$ elements. The Brandt semigroup $B(G,I)$ admits a finite basis of inverse identities but its plain identities have no finite basis.\\[.25ex]

\emph{Proof.} Since the group $\mathbb{S}_3$ is metabelian, so is $G=\mathbb{Z}\times\mathbb{S}_3$. It is known~\cite{Co67} that the group identities of any metabelian group possess a finite basis. By~\cite[Corollary~2]{Klei77},
the inverse identities of a Brandt semigroup admit a finite basis whenever so do the group identities of its structure group. Thus, we may conclude that $B(G,I)$ has a finite basis of inverse identities.

Now consider the following series of identities:
\[
L_n:\ x^{2} y_{1} \cdots y_{n} y_{n} \cdots y_{1}=y_{1} \cdots y_{n} y_{n} \cdots y_{1} x^{2},\ n=1,2,\dotsc.
\]
We aim to show that all identities $L_n$ hold in $B(G,I)$. Due to the ``if'' part of Lemma~\ref{lem:b2vg}, it amounts to verifying that they hold in both $G$ and the 5-element Brandt semigroup $B_2$. Since the group $\mathbb{S}_3$ satisfies the identity~\eqref{eq:abelian1}, this identity, which is equivalent to $L_1$, holds in $G=\mathbb{Z}\times\mathbb{S}_3$. Now it easy to verify that $G$ satisfies the identity $L_n$ by induction on $n$. Indeed, for $n>1$ we have
\begin{align*}
x^{2} y_{1}y_2 \cdots y_{n} y_{n} \cdots y_2y_{1}&=y_1(y_1^{-1}xy_1)^{2}y_2 \cdots y_{n} y_{n} \cdots y_2y_{1} &&\\
                                                 &=y_1y_2 \cdots y_{n} y_{n} \cdots y_2(y_1^{-1}xy_1)^{2}y_1 &&\text{by the inductive}\\ 
                                                 & &&\text{assumption}\\
                                                 &=y_1y_2 \cdots y_{n} y_{n} \cdots y_2y_1^{-1}x^2y_1^2 &&\\
                                                 &=y_1y_2 \cdots y_{n} y_{n} \cdots y_2y_1^{-1}y_1^2x^2 &&\text{by using~\eqref{eq:abelian1}}\\
                                                 &=y_1y_2 \cdots y_{n} y_{n} \cdots y_2y_1x^2.&&
\end{align*}
In order to show that each of the identities $L_n$ holds in $B_2=B(E,\{1,2\})$, it suffices to observe that the values of the words $x^{2} y_{1} \cdots y_{n} y_{n} \cdots y_{1}$ and $y_{1} \cdots y_{n} y_{n} \cdots y_{1} x^{2}$ under every evaluation $\varphi\colon\{x,y_1,\dots,y_n\}\to B_2$ are equal to 0 unless $x\varphi=y_k\varphi=(1,1,1)$ or $x\varphi=y_k\varphi=(2,1,2)$ for all $k=1,\dots,n$, in which case the values of these words are equal to $(1,1,1)$ or $(2,1,2)$ respectively.

Isbell~\cite{Is70} proved that no finite set of plain semigroup identities true in the groups $\mathbb{Z}$ and $\mathbb{S}_3$  implies all identities $L_n$. Hence, the plain identities of $B(G,I)$ admit no finite basis.
\hfill$\square$\\[.25ex]

Our next result also deals with the finite basis property. It immediately follows from Theorem~\ref{thm:basis}.

\begin{cor}[{\mdseries\cite[Corollary 2]{Vo85}}]
\label{cor:finite}
If a group $G$ of finite exponent admits a finite identity basis, then so does every Brandt semigroup over $G$.
\end{cor}

In particular, since every finite group possesses a finite identity basis (\cite{OaPo64}, see also \cite[Section~5.2]{Neu67}), we conclude that the plain identities of each finite Brandt semigroup have a finite basis.

\medskip

Two algebraic structures of the same type are said to be \emph{equationally equivalent} if they satisfy the same identities. Results in~\cite{Klei77}, see also \cite[Proposition~XII.4.13]{Petrich84}, imply that the following dichotomy holds for an arbitrary inverse semigroup $S$: either

(1) $S$ is equationally equivalent to an inverse semigroup that is either a group, or a group with zero, or a Brandt semigroup and that can be chosen to be finite whenever $S$ is finite, or

(2) the inverse semigroup variety generated by $S$ contains the 6-\emph{element Brandt monoid} $B_2^1$ obtained by adjoining a ``fresh'' symbol $1$ to the 5-element Brandt semigroup $B_2$ and extending the multiplication of $B_2$ so that $1$ becomes the identity element.

If $S$ and $T$ are inverse semigroups and $S$ satisfies all inverse identities of $T$, then the same holds for the plain identities of $T$ since the latter are  special instances of the former. (In the language of varieties, this means that $S$ lies in the semigroup variety generated by $T$ whenever it belongs to the inverse semigroup variety generated by $T$.) In particular, if $S$ and $T$ are equationally equivalent as inverse semigroups, they are equationally equivalent as plain semigroups as well. In view of these observations, we see that the above dichotomy persists if one considers plain semigroup identities and varieties. Thus, if $S$ is an arbitrary inverse semigroup, then either

(1') $S$ is equationally equivalent as a plain semigroup to either a group, or a group with zero, or a Brandt semigroup, each of which can be chosen to be finite whenever $S$ is finite, or

(2') the plain semigroup variety generated by $S$ contains the 6-element Brandt monoid $B_2^1$.

This dichotomy, combined with a powerful result by Sapir~\cite{Sa87}, allows us to give the following classification of finite inverse semigroups with respect to the finite basis property.
\begin{cor}[{\mdseries\cite[Corollary 3]{Vo85}}]
\label{cor:classification}
A finite inverse semigroup $S$ admits a finite basis of plain identities if and only if the plain semigroup variety generated by $S$ excludes the monoid $B_2^1$.
\end{cor}

\emph{Proof.} The ``only if'' part follows from \cite[Corollary~6.1]{Sa87}, according to which every (not necessarily inverse) finite semigroup that generates a variety containing $B_2^1$ has no finite identity basis. For the proof of the ``if'' part, we invoke the above dichotomy that allows us to assume that $S$ is either a finite group, or a finite group with zero, or a finite Brandt semigroup. We have already mentioned that every finite group possesses a finite identity basis, and so does every  finite Brandt semigroup by Corollary~\ref{cor:finite}. The remaining case of finite groups with zero easily follows from a general result by Melnik~\cite[Theorem~4]{Me70} ensuring that if a (not necessarily finite) semigroup~$T$ has a finite identity basis, then so does the semigroup $T^0$. (See \cite[Section~3]{Vo01} for a detailed explanation of how~\cite[Theorem~4]{Me70} implies this claim.) \hfill$\square$\\[.25ex]

\emph{Remark~4.} As it has been observed by Kalicki~\cite{Ka52}, there exists an algorithm to decide, given two finite algebraic structures of the same type, whether one of them belongs to the variety generated by the other. Hence, Corollary~\ref{cor:classification} provides an algorithm to decide whether or not a given finite \textbf{inverse} semigroup admits a finite basis of \textbf{plain} identities. Recall that the existence of such an algorithm remains  open for each of the following two situations: when one wants to decide whether or not a given finite \textbf{plain} semigroup admits a finite basis of \textbf{plain} identities (see \cite[Section~2]{Vo01} for a discussion) and when one wants to decide whether or not a given finite \textbf{inverse} semigroup admits a finite basis of \textbf{inverse} identities. In particular, it is not known if for a finite inverse semigroup $S$, the plain and the inverse versions of the finite basis property are equivalent. Ka\softd{}ourek~\cite{Kad03} has proved that they are equivalent provided that all subgroups of $S$ are solvable.

\section*{Acknowledgements}
The author thanks Dr.\ Ji\v{r}\'{\i} Ka\softd{}ourek who carefully examined a number of publications of the 1980s, a notable ``Sturm und Drang'' period in the theory of semigroup varieties (and corrected inaccuracies in some of these publications, see, e.g., \cite{Kad18}). In the course of his critical studies, Dr.~Ka\softd{}ourek observed a gap in~\cite{Vo85} and drew the author's attention to the fact that this gap had not been properly discussed in the literature. The present paper is a response to this fair remark.

\end{document}